
\magnification1200
\input amstex.tex
\documentstyle{amsppt}
\nopagenumbers
\hsize=12.5cm
\vsize=18cm
\hoffset=1cm
\voffset=2cm

\footline={\hss{\vbox to 2cm{\vfil\hbox{\rm\folio}}}\hss}

\def\DJ{\leavevmode\setbox0=\hbox{D}\kern0pt\rlap
{\kern.04em\raise.188\ht0\hbox{-}}D}

\def\txt#1{{\textstyle{#1}}}
\baselineskip=13pt
\def\hf{{\textstyle{1\over2}}}
\def\a{\alpha}\def\b{\beta}
\def\d{{\,\roman d}}
\def\e{\varepsilon}

\def\G{\Gamma}
\def\k{\kappa}

\def\t{\theta}
\def\={\;=\;}

\def\zt{\zeta(\hf+it)}

\def\D{\Delta}

\def\z{\zeta}

 \def\t{\theta}
\def\hf{{\textstyle{1\over2}}}
\def\txt#1{{\textstyle{#1}}}

\def\le{\leqslant} \def\ge{\geqslant}
\font\tenmsb=msbm10
\font\sevenmsb=msbm7
\font\fivemsb=msbm5
\newfam\msbfam
\textfont\msbfam=\tenmsb
\scriptfont\msbfam=\sevenmsb
\scriptscriptfont\msbfam=\fivemsb
\def\Bbb#1{{\fam\msbfam #1}}

\def \NN {\Bbb N}

\def \ZZ {\Bbb Z}

\font\ff=cmr8
\def\txt#1{{\textstyle{#1}}}
\baselineskip=13pt

\font\teneufm=eufm10
\font\seveneufm=eufm7
\font\fiveeufm=eufm5
\newfam\eufmfam
\textfont\eufmfam=\teneufm
\scriptfont\eufmfam=\seveneufm
\scriptscriptfont\eufmfam=\fiveeufm
\def\mathfrak#1{{\fam\eufmfam\relax#1}}

\font\tenmsb=msbm10
\font\sevenmsb=msbm7
\font\fivemsb=msbm5
\newfam\msbfam
     \textfont\msbfam=\tenmsb
      \scriptfont\msbfam=\sevenmsb
      \scriptscriptfont\msbfam=\fivemsb
\def\Bbb#1{{\fam\msbfam #1}}

\def \NN {\Bbb N}

\def \ZZ {\Bbb Z}

  \def\rightheadline{{\hfil{\ff
  On some results for $|\zt|$ and a divisor problem II}\hfil\tenrm\folio}}

  \def\leftheadline{{\tenrm\folio\hfil{\ff
   Aleksandar Ivi\'c and Wenguang Zhai}\hfil}}
  \def\emptyheadline{\hfil}
  \headline{\ifnum\pageno=1 \emptyheadline\else
  \ifodd\pageno \rightheadline \else \leftheadline\fi\fi}

\topmatter
\title
ON SOME MEAN VALUE RESULTS FOR THE ZETA-FUNCTION AND A DIVISOR PROBLEM II
\endtitle
\author   Aleksandar Ivi\'c and Wenguang Zhai${}^1$
 \endauthor

\nopagenumbers
\footnotetext[1]{
W. Zhai's work is  supported by the National Key Basic Research Program of China
(Grant No. 2013CB834201),
the National Natural Science Foundation of China (Grant No. 11171344),
and the Fundamental Research Funds for the
Central Universities in China (2012YS01).}

\medskip

\address
Aleksandar Ivi\'c, Katedra Matematike RGF-a
Universiteta u Beogradu, \DJ u\v sina 7, 11000 Beograd, Serbia\vskip2mm
Wenguang Zhai, Department of Mathematics,
China University of Mining and Technology,
Beijing  100083, P.R.China
\endaddress

\bigskip
\email {
\tt
aleksandar.ivic\@rgf.bg.ac.rs,  zhaiwg\@hotmail.com }\endemail
\dedicatory
\enddedicatory
\abstract
{Let $d(n)$ be the number of divisors of $n$, let
$$
\D(x) := \sum_{n\le x}d(n) - x(\log x + 2\gamma -1)
$$
denote the error term in the classical Dirichlet
divisor problem, and let $\z(s)$ denote the Riemann
zeta-function. It is shown that
$$
\int_0^T\D(t)|\zt|^2\d t \ll T(\log T)^{4}.
$$
Further,
if $2\le k\le 8$ is a fixed integer, then we prove the asymptotic formula
$$
\int_1^{T}\Delta^{k}(t)|\zeta(\hf+it)|^2\d t=c_1(k)T^{1+\frac k4}\log T+
c_2(k)T^{1+\frac k4}+O_\e(T^{1+\frac k4-\eta_k+\varepsilon}),
$$
where
$c_1(k)$ and $c_2(k)$ are explicit constants, and where
$$\eta_2=  3/20, \eta_3=  \eta_4=1/10,\ \eta_5=3/80,\ \eta_6=35/4742,\
\eta_7=17/6312,\ \eta_8=8/9433.$$
The results depend on the power moments of $\D(t)$ and $E(T)$, the classical
error term in the asymptotic formula for the mean square of $|\zt|$.
 }
\endabstract
\endtopmatter

\document

\head
1. Introduction
\endhead

As usual, let
$$
\D(x) \;:=\; \sum_{n\le x}d(n) - x(\log x + 2\gamma - 1)\qquad(x\geqslant 2)
\leqno(1.1)
$$
denote the error term in the classical Dirichlet divisor problem (see e.g.,
Chapter 3 of [4]). Also let
$$
E(T) \;:=\;\int_0^T|\zt|^2\d t - T\Bigl(\log\bigl({T\over2\pi}\bigr) + 2\gamma - 1
\Bigr)\qquad(T\geqslant2)\leqno(1.2)
$$
denote the error term in the mean square formula for $|\zt|$.
Here $d(n)$ is the number of all positive divisors of
$n$, $\z(s)$ is the Riemann zeta-function, and $ \gamma = -\G'(1) = 0.577215\ldots\,$
is Euler's constant.  In the first part of this work [9], the first author proved
several results involving the mean values of $\D(x), E(t)$ and
$$
\eqalign{
\D^*(x) :&= -\D(x)  + 2\D(2x) - \hf\D(4x)\cr&
= \hf\sum_{n\le4x}(-1)^nd(n) - x(\log x + 2\gamma - 1),\cr}
\leqno(1.3)
$$
which is the ``modified'' divisor function, introduced and studied
 by M. Jutila [12], [13].
In view of F.V. Atkinson's classical explicit formula [1] for $E(T)$,
which shows analogies between $\D(x)$ and $E(T)$, it turns out that
$\D^*(x)$ is a better analogue of $E(T)$ than $\D(x)$ itself.
Namely, M. Jutila (op. cit.) investigated both the
local and global behaviour of the difference
$$
E^*(t) \;:=\; E(t) - 2\pi\D^*\bigl({t\over2\pi}\bigr),\leqno(1.4)
$$
and in particular  he proved that
$$
\int_T^{T+H}(E^*(t))^2\d t \;\ll_\e\; HT^{1/3}\log^3T+ T^{1+\e}\qquad(1\le H\le T).
\leqno(1.5)
$$
Here and later $\e$ denotes positive constants which are arbitrarily
small, but are not necessarily the same ones at each occurrence,
while $a(x) \ll_\e b(x)$ (same as $a(x) = O_\e(b(x)))$ means that
the $|a(x)| \le Cb(x)$ for some $C = C(\e) >0, x\ge x_0$.
The significance of (1.5) is that, in view of (see e.g., [4, Chapter 15])
$$
\int_0^T(\D^*(t))^2\d t \;\sim\; AT^{3/2},\quad
 \int_0^T E^2(t)\d t \;\sim \;BT^{3/2}\quad(A,B >0,\;
T\to\infty),\leqno(1.6)
$$
it transpires that $E^*(t)$ is in the mean square sense of a lower order of magnitude than
either $\D^*(t)$ or $E(t)$. A similar mean square formula holds for $\D(t)$ as well,
and actually sharper formulas are known in all three cases; for this see
the paper of Lau--Tsang [15]. We also refer the reader
to the review paper [21] of K.-M. Tsang on this subject.

\medskip
Thus it seemed interesting to study the interplay between $\D^*(t)$ (and $\D(t)$) and $\z(s)$.
Mean values (or moments) of $|\zt|$ represent  one of the central themes in the theory of $\z(s)$,
and they have been studied extensively.
There are two monographs dedicated solely to them: the author's [7], and that of K.
Ramachandra [18].
In [9] it was proved that,
for $T^{2/3+\e} \le H = H(T) \le T$, we have
$$
 \int_T^{T+H}\D^*\bigl(\frac{t}{2\pi}\bigr)|\zt|^2\d t \;\ll\; HT^{1/6}\log^{7/2}T.
 \leqno(1.7)
$$
It was also proved that if $C$ is a suitable positive constant, then
$$
\int_0^T  {\Bigl(\D^*\bigl(\frac{t}{2\pi}\bigr)\Bigr)}^2|\zt|^2\d t =
\frac{C}{4\pi^2}T^{3/2}\Bigl(\log\frac{T}{2\pi} + 2\gamma - \frac{2}{3}\Bigr)
+ O_\e(T^{17/12+\e}),\leqno(1.8)
$$
and if $D$ is another suitable positive constant, then
$$
\int_0^T  {\Bigl(\D^*\bigl(\frac{t}{2\pi}\bigr)\Bigr)}^3|\zt|^2\d t =
DT^{7/4}\Bigl(\log\frac{T}{2\pi} + 2\gamma - \frac{4}{7}\Bigr)
+ O_\e(T^{27/16+\e}).\leqno(1.9)
$$
The proofs of (1.8) and (1.9), given in [9], exploited the special structure
of $\D^*\bigl(\frac{t}{2\pi}\bigr)$ and could not be readily extended to deal
with $\D^*(\a t)$ or $\D(\a t)$ for a given $\a >0$.

\medskip
\head
2. Statement of results
\endhead
\medskip

This paper is a continuation of the first author's paper [9] and the
second author's papers [24], [25],
where he investigated the high-power moments of $\D(x)$ and $E(t)$.

\medskip
Namely it is conjectured that the asymptotic formula
$$
\int_0^T  \D^k(t)\d t =
C_kT^{1+k/4}
+ O_\e(T^{1+k/4-c(k)+\e})\leqno(2.1)
$$
holds with an explicit constant $C_k$ and some $c(k)>0$, when $k>1$ is a given
natural number. An asymptotic formula analogous to (2.1) is also conjectured
for the moments of $E(t)$. The case $k=2$ (the mean square) of (2.1) is classic,
and it is now known that (see Lau--Tsang [15])
$$
\int_0^T  \D^2(t)\d t = C_2T^{3/2} + O(T\log^3T\log\log T),\leqno(2.2)
$$
with
$$
C_2 = (6\pi^2)^{-1}\sum_{n=1}^\infty d^2(n)n^{-3/2} = (6\pi^2)^{-1}\z^4(3/2)/\z(3)
= 0.25045\ldots,
$$
and a formula analogous to (2.2) holds for the mean square of $E(t)$.
A detailed discussion concerning the integral in (2.1) in the general case
is to be found
in the second author's paper [24], Part II, where (2.1) is established
for $5\le k \le9$, with explicit values of $c(k)$. For $k=3$
the best known value is $c(3) = 7/20$ (Ivi\'c--Sargos [11])
and for $k=4$
it is $c(4) = 3/28$ (W. Zhai [24]), and K.-L. Kong [14] has just
obtained  $c(4) = 1/8$. Like in the problem of the evaluation of
the moments $\int_0^T|\zt|^k\d t$ and similar problems, the problem
becomes progressively more difficult as $k$ increases. It is curious
that, for $2 \le k \le 9$, when it is as present known that the
asymptotic formula (2.1) holds, all the constants $C_k$ are positive for odd $k$,
implying that the values of $\D(t)$ are more biased towards
positive values. Whether this phenomenon will also happen for odd $k>9$, should
(2.1) continue to hold, is unclear.

\medskip
In this paper we are interested in a similar, but more involved problem,
namely the asymptotic evaluation of the integrals of $\D^k(t)|\zt|^2$ when
$k\in \NN$ is fixed. We succeeded in applying the existing results on the
moments of $\D(t)$ and $E(t)$ to the evaluation of the integrals of $\D^k(t)|\zt|^2$.
Our methods at present work for $1\le k\le 8$, and
the results are as follows.

\medskip
THEOREM 1. {\it We have}
$$
\int_0^T\D(t)|\zt|^2\d t \;\ll\; T(\log T)^{4}.\leqno(2.3)
$$
\medskip
THEOREM 2. {\it If $k$ is a fixed integer for which $2\le k \le 8$, then we have}
$$
\int_1^{T}\Delta^{k}(t)|\zeta(\hf+it)|^2\d t=c_1(k)T^{1+\frac k4}\log T+
c_2(k)T^{1+\frac k4}+O_\e(T^{1+\frac k4-\eta_k+\varepsilon}),\leqno(2.4)
$$
{\it where
$c_1(k)$ and $c_2(k)$ are explicit constants, and where}
$$
\eta_2=  \eta_3=  \eta_4=1/10,\; \eta_5=3/80,\; \eta_6=35/4742,\;
 \eta_7=17/6312,\; \eta_8=8/9433.
$$

\medskip
Note that the values of $\eta_2,  \eta_3,  \eta_4$ in Theorem 2 are identical,
which is due to the general argument used in the proof in Section 5. However,
we can combine the arguments of Theorem 1 and Theorem 2 to obtain
improvements on the values of $\eta_2$ and $\eta_3$. We shall give
the details only for $\eta_2$, while the case of $\eta_3$ is technically
quite complicated. The result is

\medskip
THEOREM 3. {\it When $k=2$, we can take $\eta_2= 3/20$ in Theorem} 2.

\medskip
{\bf Remark 1}. Note that, for $H=T$, (2.3) improves (1.7) a lot. It is an open
problem to find the lower bound for the integral in (2.3), since it is well-known
that $\D(x)$ changes sign in every interval of the form $[T, T + A\sqrt{T}]$
for a suitable $A>0$ and $T\ge T_0$
(see the first author's paper [6]). On the other hand, one has (by (3.1) and
(3.7) of [4]) that
$$
\int_1^X\D(x)\d x= \txt{1\over4}X + O(X^{3/4}).
$$
Using this formula it may be conjectured that
$$
\int_1^T\D(t)|\zt|^2\d t =
{T\over4}\Bigl(\log\frac{T}{2\pi} + 2\gamma - 1\Bigr)
+ O_\e(T^{3/4+\e}),\leqno(2.5)
$$
however obtaining any asymptotic formula for the integral in (2.5)
is difficult.

\medskip
{\bf Corollary 1}. We also have
$$
 \int_0^TE^*(t)|\zt|^2\d t \;\ll\; T(\log T)^4.\leqno(2.6)
 $$
This follows from (1.4), (2.1) (since it will hold with $\D^*(t/(2\pi))$
instead od $\D(t)$),
$$
\int_0^TE(t)|\zt|^2\d t = \pi T\Bigl(\log\frac{T}{2\pi} + 2\gamma-1\Bigr) + U(T),
\leqno(2.7)
$$
where
$$
U(T) = O(T^{3/4}\log T),\quad U(T) = \Omega_\pm(T^{3/4}\log T).\leqno(2.8)
$$
\medskip
The asymptotic formulas (2.7)--(2.8) are due to the first author [5]. They show,
up to the numerical constants which are involved, the
true order of magnitude of the function $U(T)$. Here the symbol $f(x) = \Omega_\pm(g(x))$
has its standard meaning, namely that both $\,\limsup\limits_{x\to\infty}f(x)/g(x) > 0\,$ and
$\,\liminf\limits_{x\to\infty}f(x)/g(x)<0\,$ holds.
\medskip
The analogy between (2.5) and (2.7) is obvious, however the latter is much less
difficult. Namely the defining relation (1.2) yields, by differentiation,
$$
|\zt|^2 = \log\bigl(\frac{t}{2\pi}\bigr) + 2\gamma + E'(t),\leqno(2.10)
$$
and one can easily integrate $E^k(t)E'(t)\;(k\in\NN)$. Thus the integral in (2.4) is
more difficult to evaluate than the corresponding problem when $\D^k(t)$ is
replaced by $E^k(t)$.

{\bf Remark 2}.
The methods of proof of (2.4) allow one
to carry over the results of Theorem 1, Theorem 2 and Theorem 3 to
integrals where $\D(t)$ is replaced by $\D(\a t)$ or $\D^*(\a t)$ for
any given $\a>0$.

\medskip
{\bf Remark 3}.
It would be interesting to analyze the error term in (2.4) and see how small it
can be, i.e., to obtain an omega-result (recall that $f(x) = \Omega(g(x))$ means
that $f(x) = o(g(x))$ does not hold as $x\to\infty$).

\medskip
{\bf Remark 4}.  For $k=2$ one can compare (2.4) with the
corresponding result of the first author [5],
where it was obtained that
$$
\int_0^T E ^2(t)|\zt|^2\d t =
D_2T^{3/2}\Bigl(\log\frac{T}{2\pi} + 2\gamma - \frac{2}{3}\Bigr)
+ O(T\log^6T),
$$
where
$$
D_2 \;=\; \frac{2\z^4(3/2)}{3\sqrt{2\pi}\z(3)}.
$$

\medskip
{\bf Remark 5}. Finally we indicate two possible generalizations
of our results. Let, as usual, $r(n) = \sum_{n=a^2+b^2}1$ denote
the number of ways $n$ may be represented as a sum of two integer squares, and
let  $\varphi(z)$ be a
holomorphic cusp form of weight $\kappa$ with respect to the full
modular group $SL(2,\ZZ)$, and denote by $a(n)$ the $n$-th Fourier
coefficient of $\varphi(z)$. We suppose that $\varphi(z)$ is a
normalized eigenfunction for the Hecke operators $T(n)$, that is,
$  a(1)=1  $ and $  T(n)\varphi=a(n)\varphi $ for every $n \in
\NN$.  The classical example is $a(n) = \tau(n)\;(\k=12)$,
 the Ramanujan $\tau$-function defined by
$$
\sum_{n=1}^\infty \tau(n)x^n \=
x{\left\{(1-x)(1-x^2)(1-x^3)\cdots\right\}}^{24}\qquad(\,|x| < 1).
$$
If $P(x) := \sum\limits_{n\le x }r(n) - \pi x$ denotes then the error term in
the classical circle problem and $A(x) := \sum\limits_{n\le x}a(n)$, then Theorem 2
can be generalized to integrals
$$
\int_0^T P^k(t)|\zt|^2\d t,\quad \int_0^T A^k(t)|\zt|^2\d t,\leqno(2.11)
$$
more precisely if $A(t)$ replaced by the normalized
function $A^*(t) := \sum\limits_{n\le t}a(n)n^{{1-\k\over2}}$,
since $a(n)$ behaves similarly to $n^{(\kappa-1)/2}d(n)$. For the analogues of
Lemma 4 to $P(x)$ and $A(x)$ the reader should see e.g., section 3 of [8].
The analogue of (3.1) for $\D(x)$ will hold with a poorer $\t$ (with $\t = 1/3$
in case of $A^*(x)$), and the analogues of the exponents $\eta_k$ will not be
as good as those of Theorem 2.

\medskip

\head
3. The necessary lemmas
\endhead
In this section we shall state some lemmas needed for the proof of our theorems.
The proofs of the theorems themselves will be given in Section 4, Section 5
 and Section 6.
\medskip

\medskip
LEMMA 1. {\it
There exists a constant $\t$ such that $1/4 \le \theta<1/3$ and
$$\Delta(x)\ll_\e x^{\theta+\varepsilon}, \ \ E(t)\ll_\e
t^{\theta+\varepsilon}.\leqno(3.1)
$$
In particular,  we can take $\theta=131/416=0.3149\cdots.$}

\medskip
The proofs of the bounds in (3.1) are due to M.N. Huxley [3] and N. Watt [23], respectively,
and they are the sharpest ones known.
It is commonly conjectured that $\t = 1/4$ is permissible, but this is out of reach
at present. It is known that $\t < 1/4$ cannot hold (see e.g., [4], Chapter 13 and Chapter 15).

\medskip
LEMMA 2.
{\it Suppose $\theta$ is the constant in Lemma 1. Then for any $A$
satisfying $0\le A\le 11$ we have
$$\int_1^T |\Delta(x)|^A\d x\ll_\e T^{1+M(A)+\varepsilon}\leqno(3.2)
$$
and
$$\int_1^T |E(t)|^A\d t\ll_\e T^{1+M(A)+\varepsilon},\leqno(3.3)$$
where
$$M(A):=\max\left(\frac{A}{4}, \theta(A-2)\right).\leqno(3.4)$$}

\medskip
We note that, for real $k \in [0,9]$, the limits
$$
E_k
\;:=\; \lim_{T\to\infty}T^{-1-k/4}\int_0^T|E(t)|^k\d t\leqno(3.5)
$$
exist. The analogous result holds also for the moments of $\D(t)$.
This was proved by
D.R. Heath-Brown [2], who used (3.4) in his proof. He also showed that
the limits of moments (both of $\D(t)$ and $E(t)$) without
absolute values also exist when  $k=1,3,5,7$ or 9.
The merit of (3.5) that it gets rid of ``$\e$'' and establishes the existence
of the limit (but without an error term).

\medskip
LEMMA 3. {\it We have
$$
\int_0^T|\zt|^4 \d t = TQ_4(\log T) + O(T^{2/3}\log^8T),\leqno(3.6)
$$
where $Q_4(x)$ is an explicit polynomial of degree four in $x$ with
leading coefficient $1/(2\pi^2)$.
}

\medskip
This result was proved first (with error term $O(T^{2/3}\log^CT)$) by
Y. Motohashi and the author [10]. The value $C=8$ was given by Y. Motohashi
in his monograph [17]. We shall not need the full strength of (3.6), but just the
upper bound $O(T\log^4T)$ for the integral in question.

\medskip
LEMMA 4. {\it For $1 \le N \ll x$ we have
$$
\D(x) = {1\over\pi\sqrt{2}}x^{1\over4}
\sum_{n\le N}d(n)n^{-{3\over4}}
\cos(4\pi\sqrt{nx} - {\txt{1\over4}}\pi) +
O_\e(x^{{1\over2}+\e}N^{-{1\over2}}).
\leqno(3.7)
$$
and}
$$
\D^*(x) = {1\over\pi\sqrt{2}}x^{1\over4}
\sum_{n\le N}(-1)^nd(n)n^{-{3\over4}}
\cos(4\pi\sqrt{nx} - {\txt{1\over4}}\pi) +
O_\e(x^{{1\over2}+\e}N^{-{1\over2}}).
\leqno(3.8)
$$
\medskip
The expression (3.8) for $\D^*(x)$ (see [4], Chapter 15) is the analogue
of the classical truncated Vorono{\"\i} formula (3.7) for $\D(x)$ (ibid. Chapter 3),
only the sum in the expression for $\D^*(x)$ has an additional factor $(-1)^n$.
Actually G.F. Vorono{\"\i} [22] proved long ago an explicit formula for $\D(x)$ as
a series containing
the Bessel functions $K_1$ and $Y_1$
(see e.g., [4], Chapter 3). However, to avoid the questions
of convergence it is in practice usually more expedient to work with (3.7), which
is sufficient for many purposes.

\medskip

LEMMA 5. {\it For $Q \gg x \gg 1$ we have
$$
\D(x) = {1\over\pi\sqrt{2}}x^{1\over4}
\sum_{n\le Q}d(n)n^{-{3\over4}}
\cos(4\pi\sqrt{nx} - {\txt{1\over4}}\pi) + F(x),\leqno(3.9)
$$
where
$F(x) \ll x^{-1/4}$ if $||x||\gg x^{5/2}Q^{-1/2}$, and  we always
have $F(x) \ll_\e x^\e$.}

\medskip
This result ($||x||$ denotes as usual the distance of $x$ to the
nearest integer) is due to T. Meurman [16]. It shows that, unless $x$
is close to an integer, the error term in the
 truncated Vorono{\"\i} formula for $\D(x)$  is small.

\medskip
LEMMA 6.  {\it Let $k\ge 2$ be a fixed
integer and $\delta > 0$ be given.
Then the number of integers $n_1,n_2,n_3,n_4$ such that
$N < n_1,n_2,n_3,n_4 \le 2N$ and}
$$
|n_1^{1/k} + n_2^{1/k} - n_3^{1/k} - n_4^{1/k}| \;<\; \delta N^{1/k}
$$
{\it is, for any given $\e>0$,}
$$
\ll_\e\; N^\e(N^4\delta + N^2).\leqno(3.10)
$$

\medskip
Lemma 6 was proved by analytic methods
by  Robert--Sargos [19]. When $k=2$, it represents a powerful arithmetic
tools which is essential in the analysis when the
biquadrate of exponential sums involving $\sqrt n$ appears.

\medskip
LEMMA 7. {\it We have
$$
\sum_{n\le x}d^2(n) \;=\; \frac{1}{\pi^2}x\log^3x + O(x\log^2x).\leqno(3.11)
$$
}
\medskip
This is a well-known elementary formula; see e.g., page 141 of [4]. It follows from the
series representation
$$
\sum_{n=1}^\infty d^2(n)n^{-s} = \frac{\z^4(s)}{\z(2s)}\qquad(\Re s > 1).
$$

\medskip
LEMMA 8. {\it For $1\le r \ll x$ we have}
$$
\sum_{n\le x}d(n)d(n+r) \ll \sum_{d|r}\frac1d \cdot x\log^2x.
$$

\medskip
This follows e.g., from a theorem of P. Shiu [20] on multiplicative functions.

\medskip
LEMMA 9 . {\it Let $0 < A < A'$ be any two fixed constants
such that $AT < N < A'T$, and let $N' = N'(T) =
T/(2\pi) + N/2 - (N^2/4+ NT/(2\pi))^{1/2}$. Then}
$$
E(T) = \Sigma_1(T) + \Sigma_2(T) + O(\log^2T),
$$
{\it where}
$$
\Sigma_1(T) = 2^{1/2}(T/(2\pi))^{1/4}\sum_{n\leqslant N}(-1)^nd(n)n^{-3/4}
e(T,n)\cos(f(T,n)),
$$
$$
\Sigma_2(T) = -2\sum_{n\leqslant N'}d(n)n^{-1/2}{\Bigl(\log {T\over2\pi n}\Bigr)}^{-1}
\cos\left(T\log \Bigl( {T\over2\pi n}\Bigr) - T + {\txt{1\over4}}\pi\right),
$$
{\it with}
$$
\eqalign{\cr&
f(T,n) = 2T{\roman {arsinh}}\,\Bigl(\sqrt{\pi n\over2T}\,\Bigr) + \sqrt{2\pi nT
+ \pi^2n^2} - {\txt{1\over4}}\pi\cr&
=  -\txt{1\over4}\pi + 2\sqrt{2\pi nT} +
\txt{1\over6}\sqrt{2\pi^3}n^{3/2}T^{-1/2} + a_5n^{5/2}T^{-3/2} +
a_7n^{7/2}T^{-5/2} + \ldots\,,\cr}
$$
$$\eqalign{\cr
e(T,n) &= (1+\pi n/(2T))^{-1/4}{\Bigl\{(2T/\pi n)^{1/2}
{\roman {arsinh}}\,\Bigl(\sqrt{\pi n\over2T}\,\Bigr)\Bigr\}}^{-1}\cr&
= 1 + O(n/T)\qquad(1 \leqslant n < T),
\cr}
$$
{\it and $\,{\roman{ar\,sinh}}\,x = \log(x + \sqrt{1+x^2}\,).$}

\medskip
This is the famous formula of F.V. Atkinson [1]; proofs can be also
found in [4] and [7].

\medskip
LEMMA 10. {\it Let $p_1, p_2, \ldots ,p_r >0$ and $f_1(x), f_2(x), \ldots, f_r(x)\ge 0$
be continuous functions in $[a,b]\;(a<b)$. Then if
$$
\frac{1}{p_1} + \frac{1}{p_2} + \ldots + \frac{1}{p_r} =1,
$$
we have
$$
\int_a^b f_1(x)\ldots f_r(x)\d x \le {\left(\int_a^b {f}^{p_1}_1(x)\d x\right)}^{1/p_1}\cdots
{\left(\int_a^b {f}^{p_r}_r(x)\d x\right)}^{1/p_r}.
$$
}

\medskip
This is the classical H\"older inequality for integrals, the case $r=2, p_1 = p_2 = 1/2$ being the
Cauchy-Schwarz inequality. It will be repeatedly used in the proofs.
\medskip
\head
4. Proof of  Theorem 1
\endhead

\medskip
It suffices to consider in (2.3) the integral over $[T,\,2T]$,
to replace then $T$ by $T2^{-j}\;(j=1,2,\ldots)$ and sum the
resulting estimates.
We suppose $T\le t \le 2T$, take $Q = T^7$ in  Lemma 5 and write
$$
\D(t) =
\D_1(t)+ \D_2(t) + F(t),\leqno(4.1)
$$
where $F(t)$ is as in Lemma 5, and
$$
\eqalign{
\D_1(t) &:= \frac{t^{1/4}}{\pi\sqrt{2}}\sum_{n\le T}
d(n)n^{-3/4}\cos(4\pi\sqrt{nt}-\textstyle{\frac \pi 4}),\cr
\D_2(t) &:= \frac{t^{1/4}}{\pi\sqrt{2}}\sum_{T<n\le Q}
d(n)n^{-3/4}\cos(4\pi\sqrt{nt}-\textstyle{\frac \pi 4}).\cr}
\leqno(4.2)
$$
Therefore
$$
\int_T^{2T}\D(t)|\zt|^2\d t = \int_T^{2T}(\D_1(t)+\D_2(t) + F(t))|\zt|^2\d t.
$$

By the Cauchy-Schwarz inequality (Lemma 10) it is seen that
the term $F(t)$ in (4.1) makes a contribution of $O(T^{3/4}\log T)$.
The contribution containing $\D_2(t)$ is, by the first derivative test
(see e.g., Lemma 2.1 of [4]), Lemma 3 and Lemma 7,
$$
\eqalign{
&
\ll T^{1/4}{\Bigl\{\int_T^{2T}{\Bigl|\sum\limits_{T<n\le Q}
d(n)n^{-3/4}\cos(4\pi\sqrt{nt}-\textstyle{\frac \pi 4})\Bigr|}^2\d t
\int_T^{2T}|\zt|^4\d t\Bigr\}}^{\frac12}\cr&
\ll T^{3/4}\log^2T{\Biggl\{T^{1/2}\log^3T
+ T^{1/2}\sum_{T<m\ne n\le Q}\frac{d(m)d(n)}{(mn)^{3/4}|\sqrt{m}-\sqrt{n}|}
\Biggr\}}^{\frac12}.\cr}\leqno(4.3)
$$
In the double sum in in (4.3), the contribution when $m\ge4n$ 
or $n\ge4m$ is $\ll \log^3T$.
The contribution of the remaining terms is, supposing $m>n$, setting $m=n+r$
and using Lemma 8,
$$
\ll \sum_{r\ll Q}\frac1r \sum_{n\le Q}\frac{d(n)d(n+r)}{n} \ll \sum_{r\ll Q}
\frac1r\sum_{d|r}\frac1d \log^3T\ll \log^4T.
$$
Therefore the contribution containing $\D_2(t)$ is
$$
\ll T^{3/4}\log^2T(T^{1/2}\log^4T)^{1/2} = T\log^4T.
$$
Further we have, by (2.10),
$$
\eqalign{
&\int_T^{2T}\D_1(t)|\zt|^2\d t\cr&
= \int_T^{2T}\frac{t^{1/4}}{\pi\sqrt{2}}\sum_{n\le T}
d(n)n^{-3/4}\cos(4\pi\sqrt{nt}-\textstyle{\frac \pi 4})
\Bigl(\log\frac{t}{2\pi}+2\gamma + E'(t)\Bigr)\d t\cr&
= I_1(T) + I_2(T),
\cr}\leqno(4.4)
$$
say, where by the first derivative test
$$
\eqalign{
I_1(T):&=\int_T^{2T}\frac{t^{1/4}}{\pi\sqrt{2}}\sum_{n\le T}
d(n)n^{-3/4}\cos(4\pi\sqrt{nt}-\textstyle{\frac \pi 4})
\Bigl(\log\frac{t}{2\pi}+2\gamma \Bigr)\d t
\cr&
\ll T^{1/4}\log T\cdot \sum_{n\le T}d(n)n^{-3/4}T^{1/2}n^{-1/2}\cr&
\ll T^{3/4}\log T,
\cr}
$$
since $\sum\limits_{n\ge1}d(n)n^{-\a}$ converges for $\a>1$. The integral $I_2(T)$, namely
$$
I_2(T) := \int_T^{2T}E'(t)\frac{t^{1/4}}{\pi\sqrt{2}}\sum_{n\le T}
d(n)n^{-3/4}\cos(4\pi\sqrt{nt}-\textstyle{\frac \pi 4})\d t,
$$
is integrated by parts. Since $E(t) \ll t^{1/3}$ (see e.g., Chapter 15 of [4], also
follows trivially from Lemma 1), the integrated terms are trivially
$$
\ll\; T^{\frac{1}{3}+\frac{1}{4}}T^{\frac{1}{4}}\log T \;\ll\; T^{\frac{5}{6}}\log T.
$$
There remains a multiple of
$$
\eqalign{
&-\frac{1}{4}\int_T^{2T}t^{-3/4}E(t)\sum_{n\le T}
d(n)n^{-3/4}\cos(4\pi\sqrt{nt}-\textstyle{\frac \pi 4})\d t\cr&
+2\pi \int_T^{2T}t^{-1/4}E(t)\sum_{n\le T}
d(n)n^{-1/4}\sin(4\pi\sqrt{nt}-\textstyle{\frac \pi 4})\d t.
\cr}\leqno(4.5)
$$
Both integrals in (4.5) are estimated analogously, and clearly it is the
latter which is larger. We replace $E(t)$ by the expression given by
Atkinson's formula (see Lemma 9). Thus, taking $N=T$ in Atkinson's formula,
$$
\int_T^{2T}t^{-1/4}E(t)\sum_{n\le T}
d(n)n^{-1/4}\sin(4\pi\sqrt{nt}-\textstyle{\frac \pi 4})\d t = J_1(T) + J_2(T) + J_3(T),
$$
say, where
$$
\eqalign{
J_1(T) :&= \int_T^{2T}t^{-1/4}\sum\nolimits_1(t)\sum_{n\le T}
d(n)n^{-1/4}\sin(4\pi\sqrt{nt}-\textstyle{\frac \pi 4})\d t,\cr
J_2(T) :&= \int_T^{2T}t^{-1/4}\sum\nolimits_2(t)\sum_{n\le T}
d(n)n^{-1/4}\sin(4\pi\sqrt{nt}-\textstyle{\frac \pi 4})\d t,\cr
J_3(T) :&= \int_T^{2T}t^{-1/4}O(\log^2T)\sum_{n\le T}
d(n)n^{-1/4}\sin(4\pi\sqrt{nt}-\textstyle{\frac \pi 4})\d t.
\cr}
$$
The Cauchy-Schwarz inequality gives
$$
\eqalign{&
J_2(T) + J_3(T) \cr&
\ll T^{-1/4}\left\{\int_T^{2T}
\left(\sum\nolimits_2^2(t)+\log^4T\right)\d t\int_T^{2T}
\Bigl|\sum_{n\le T}d(n)n^{-1/4}e^{4\pi i\sqrt{nt}\,}\Bigr|^2\d t\right\}^{1/2}.
\cr}\leqno(4.6)
$$
But (see Chapter 15 of [4])
$$
\int_T^{2T}\sum\nolimits_2^2(t)\d t \;\ll\;T\log^4T,
$$
since $\sum_2(t)$ is essentially a Dirichlet polynomial of length $\asymp T$.
In the other integral in (4.6) we square out the sum and integrate. The contribution
is
$$
\eqalign{&
\ll T\sum_{n\le T}d^2(n)n^{-1/2}+ \sum_{m\ne n\le T}\frac{d(m)d(n)}{(mn)^{1/4}}
\int_T^{2T}e^{4\pi i(\sqrt{m}-\sqrt{n})\sqrt{t}}\d t\cr&
\ll T^{3/2}\log^3T + T^{1/2}\sum_{m\ne n\le T}\frac{d(m)d(n)}{(mn)^{1/4}
|\sqrt{m}-\sqrt{n}|},
\cr}
$$
by the first derivative test and Lemma 7. Note that
if $m\le n/2$, then $|\sqrt{m}-\sqrt{n}|^{-1} \ll n^{-1/2}$, while if $m > 2n$,
then $|\sqrt{m}-\sqrt{n}|^{-1} \ll m^{-1/2}$. When $m \asymp n$ the contribution is
estimated, as in (4.3), by Lemma 6. In this way it is seen that
$$
\int_T^{2T}
\Bigl|\sum_{n\le T}d(n)n^{-1/4}e^{4\pi i\sqrt{nt}}\Bigr|^2\d t \;\ll\; T^{3/2}\log^3T,
\leqno(4.7)
$$
and we obtain
$$
J_2(T) + J_3(T) \ll T^{-1/4}T^{1/2}\log^2T\cdot T^{3/4}\log^{3/2}T \ll T(\log T)^{7/2}.
$$
It remains to deal with ($c$ is a constant)
$$
J_1(T) = c\int_T^{2T}\sum_{m\le T}\frac{(-1)^md(m)}{m^{3/4}}
e(t,m)\cos f(t,m)\sum_{n\le T}\frac{d(n)}{n^{1/4}}
\cos(4\pi\sqrt{nt}-\textstyle{\frac \pi 4})\d t.
$$
We split the sums over $m,n$ into $O(\log^2T)$ subsums with the ranges of summation
$M<m\le M'\le2M, N<n\le N'\le2N$, respectively.
We write the cosines as exponentials and then
obtain $\ll \log^2T$ sums of the form
$$
\eqalign{&
\sum_{M< m\le M'}\frac{(-1)^md(m)}{m^{3/4}}\sum_{N<n\le N'}\frac{d(n)}{n^{1/4}}\times
\cr&
\times
\int_T^{2T}e(t,m)\exp
\Bigl(4\pi i\sqrt{nt} - i\sqrt{8\pi mt}-ia_3m^{3/2}t^{-1/2}-\ldots\Bigr)\d t.\cr}\leqno(4.8)
$$
There is also the expression with $+$ in place of $-$ in the exponential,
and their conjugates, but
it is (4.8) that is the relevant sum. The smooth function
$e(t,m) \,(= 1 + O(m/T))$ may be removed
on applying integration by parts. Furthermore, if $N \ge 100M$, then by the first
derivative test the contribution of the expression in (4.8) is
$$
\ll T^{1/2}\sum_{M<m\le2M}d(m)m^{-3/4}\sum_{N<n\le2N}d(n)n^{-3/4}
\ll T\log^2T,\leqno(4.9)
$$
and the same bound as in (4.9) holds when $M \ge 100N$. These sums in total make
a contribution which is $\ll T\log^4T$.

\medskip
There remains the case when $N/100 < M < 100N$. Then we use the Cauchy-Schwarz inequality
for integrals. The contribution is
$$
\eqalign{
\ll \;&\Biggl
\{\int_T^{2T}\Bigl|\sum_{N<n\le N'}\frac{d(n)}{n^{1/4}}e^{4\pi i\sqrt{nt}\,}\Bigr|^2\d t\cr&
\times
 \int_T^{2T}\Bigl|\sum_{M<m\le M'}\frac{(-1)^md(m)}{m^{3/4}}e(t,m)e^{if(t,m)}\Bigr|^2\d t\Biggr\}^{1/2}.
 \cr}\leqno(4.10)
$$
Here the first integral is estimated as in (4.7), more precisely by
$$
O(TN^{1/2}\log^3T + T^{1/2}N\log^3T).
$$
The second integral  is, by the first derivative test and Lemma 7,
$$
\ll T\sum_{m\ge M}\frac{d^2(m)}{m^{3/2}} +
\sum_{M<k\ne m\le M'}\frac{d(k)d(m)e(t,k)e(t,m)}{(km)^{3/4}}\max_{ t\in[T,2T]}\frac{1}{|f'(t,m)-f'(t,k)|}.
$$
We have
$$
f'(t,\ell) = \frac{\partial f(t,\ell)}{\partial t} = 2\roman{ar\,sinh}\sqrt{\frac{\pi \ell}{2t}},
$$
so that by the mean value theorem we obtain
$$
|f'(t,m)-f'(t,k)| \asymp \frac{|\sqrt{k}-\sqrt{m}|}{\sqrt{T}}\qquad(k\ne m,\; T\le t \le2T).
$$
Hence the last expression above is
$$
\ll TM^{-1/2}\log^3T + T^{1/2}\log^4T.
$$
It is seen then, since $M\asymp N$, that the expression in (4.10) is
$$
\eqalign{&
\ll{\Bigl((TM^{1/2}\log^3T + T^{1/2}M\log^3T)(TM^{-1/2}\log^3T +
T^{1/2}\log^4T)\Bigr)}^{1/2}\cr&
\ll (T^2\log^6T + T^{3/2}M^{1/2}\log^7T)^{1/2}.\cr}
$$
Taking $M = T2^{-j}$ and summing over $j$ we obtain that the
contribution of $J_1(T)$ is $O(T\log^{4}T)$, since $M\asymp N$ in
the relevant cases. This gives
$$
\int_T^{2T}\D(t)|\zt|^2\d t \ll T\log^4 T,
$$
and thus completes the proof of Theorem 1.

\medskip
\head
5. The proof of Theorem 2
\endhead

\medskip

Like in the proof of Theorem 1 it suffices to prove the result
for the integral over $[T,\,2T]$, where $T \,(\ge 10)$ is large.
Henceforth let
$$
T \le t \le 2T,\quad T^{1/2} \ll y = y(T) \ll T,\leqno(5.1)
$$
where $y$ will be determined later.
Write
$$
\Delta(t)=\Delta_1(t,y)+\Delta_2(t,y),\leqno(5.2)
$$
where
$$
\Delta_1(t,y):=\frac{t^{1/4}}{\sqrt 2\pi}
\sum_{n\le y}\frac{d(n)}{n^{3/4}}\cos(4\pi\sqrt{nt}-\textstyle{\frac{\pi}{4}}),
\leqno(5.3)
$$
and by Lemma 4 (with $N =y$)
$$
\D_2(t,y) \ll_\e T^{1/2+\e}y^{-1/2}\quad (\ll_\e T^{1/4+\e}).\leqno(5.4)
$$
Then we have
$$
\eqalign{
&\int_T^{2T}\Delta^{k}(t)|\zt|^2\d t\cr&
=\int_T^{2T}{\Bigl(\Delta_1(t,y)+\Delta_2(t,y)\Bigr)}^k|\zt|^2\d t\cr
&=\int_1\;+\;O\left(\int_2+\int_3\right),\cr}\leqno(5.5)
$$
where
$$
\eqalign{
&\int_1:=\int_T^{2T}\Delta_1^{k}(t,y)|\zt|^2\d t,\cr
&\int_2:=\int_T^{2T}|\Delta_1^{k-1}(t,y)\Delta_2(t,y)||\zt|^2\d t,\cr
&\int_3:=\int_T^{2T}|\Delta_2(t,y)|^{k}|\zt|^2\d t.\cr}
$$
In order to estimate $\int_3$ we need (5.4) and
$$
\int_T^{2T}|\Delta_2(t,y)|^2\d t \ll_\e T^{3/2+\varepsilon}y^{-1/2},\leqno(5.6)
$$
which follows as in the proof of (4.3). From (5.4), (5.6),  the fourth power moment
of $\zt$ and the Cauchy-Schwarz inequality (Lemma 10) we obtain
$$
\eqalign{
\int_3&\;\ll_\e\; \left(\frac{T^{1/2+\varepsilon}}{y^{1/2}}\right)^{k-1}\int_T^{2T}
|\Delta_2(t,y)| |\zt|^2\d t\cr
&\;\ll_\e\; \left(\frac{T^{1/2+\varepsilon}}{y^{1/2}}\right)^{k-1}
\left(\int_T^{2T}|\Delta_2(t,y)|^2\d t  \right)^{1/2}\left(\int_T^{2T}|
\zt|^4\d t  \right)^{1/2}\cr
&\;\ll_\e\; T^{\frac{2k+3}{4}+\varepsilon}y^{-\frac{2k-1}{4}}.\cr}\leqno(5.7)
$$

Now we evaluate $\int_1.$ We write (2.10) as
$$
|\zt|^2=\log t+C+E^{\prime}(t),
$$
where henceforth we set $C=2\gamma-\log 2\pi$  for brevity.
Therefore we have
$$
\eqalign{
\int_1&=\int_T^{2T}\Delta_1^{k}(t,y)(\log t+C)\d t+\int_T^{2T}
\Delta_1^{k}(t,y)E^{\prime}(t) \d t\cr
&=\int_{11}+\int_{12},\cr}\leqno(5.8)
$$
say. We  bound first $\int_{12}.$ Using  integration by parts and Lemma 1 we obtain
$$
\eqalign{
\int_{12}&=\Delta_1^{k}(t,y)E(t)\Bigl|_T^{2T}-
k \int_T^{2T}\Delta_1^{k-1}(t,y)\Delta_1^{\prime}(t,y)E(t) \d t\cr
&\ll_\e T^{(k+1)\theta+\varepsilon}+  \left|\int_{12}^{*} \right|,\cr}
\leqno(5.9)
$$
say, where
$$
\int_{12}^{*}:=\int_T^{2T}\Delta_1^{k-1}(t,y)\Delta_1^{\prime}(t,y)E(t)\d t.
$$

In order to bound $\int_{12}^{*}$, we need upper bounds for
the second and the fourth moment of $\Delta_1^{\prime}(t,y).$ It is easily seen that
$$
\eqalign{
 \Delta_1^{\prime}(t,y)&=\frac{t^{-3/4}}{4\sqrt 2\pi}\sum_{n\le y}\frac{d(n)}{n^{3/4}}
 \cos(4\pi\sqrt{nt}-\textstyle{\frac{\pi}{4}})\cr
 & -\sqrt 2 t^{-1/4}\sum_{n\le y}\frac{d(n)}{n^{1/4}}\sin(4\pi\sqrt{nt}-\textstyle{\frac{\pi}{4}})
 \cr
&\ll t^{-1}|\Delta_1(t,y)|+t^{-1/4}\left|\sum_{n\le y}\frac{d(n)}{n^{1/4}}e(2\sqrt{nt} )\right|.\cr}
\leqno(5.10)
$$
Since  $\Delta_2(t,y)\ll T^{1/4+\varepsilon},$ by (5.4), it follows that
$$
\Delta_1(t,y)\ll_\e |\Delta(t)|+T^{1/4+\varepsilon}.
$$
Thus by Lemma 2 we have, for any $0\le A\le 11$,
$$
\eqalign{
\int_T^{2T}|\Delta_1(t,y)|^A\d t &\ll_\e \int_T^{2T}|\Delta(t)|^A\d t+T^{1+A/4+\varepsilon}\cr
&\ll_\e T^{1+M(A)+\varepsilon},
\cr}\leqno(5.11)
$$
where $M(A)$ is as in (3.4) of Lemma 2.
For the mean square of $\Delta_1^{\prime}(t,y)$ we have, by (5.10),
$$
\eqalign{
&\int_T^{2T}|\Delta_1^{\prime}(t,y)|^2\d t\ll \int_T^{2T}t^{-2}|\Delta_1(t,y)|^2\d t\cr
& +
\int_T^{2T}t^{-1/2}\left|\sum_{n\le y}\frac{d(n)}{n^{1/4}}e(2\sqrt{nt} )\right|^2\d t
\cr
&\ll T^{-1/2}+T^{-1/2}\sum_{m,n\le y}\frac{d(m)d(n)}
{(mn)^{1/4}}\left|\int_T^{2T}e(2(\sqrt m-\sqrt n)\sqrt t\,)\d t\right|\cr
&\ll T^{1/2}\sum_{n\le y}\frac{d^2(n)}{n^{1/2}}+T^{-1/2}\sum_{m\ne n\le y}
\frac{d(m)d(n)}{(mn)^{1/4}}\left|\int_T^{2T}e(2(\sqrt m-\sqrt n)\sqrt t\,)\d t\right|\cr
&\ll  T^{1/2}\sum_{n\leq y}\frac{d^2(n)}{n^{1/2}}+ \sum_{m\ne n\le y}
\frac{d(m)d(n)}{(mn)^{1/4}|\sqrt m-\sqrt n|}\cr
&\ll (yT)^{1/2}\log^3 T,\cr}\leqno(5.12)
$$
where we used the first derivative test and Lemma 8.

For the fourth moment of $\Delta_1^{\prime}(t,y)$ we have, by  (5.10), that
$$
\eqalign{
&\int_T^{2T}|\Delta_1^{\prime}(t,y)|^4dt\ll \int_T^{2T}t^{-4}|\Delta_1(t,y)|^4\d t\cr
& +
\int_T^{2T}t^{-1 }\left|\sum_{n\leq y}\frac{d(n)}{n^{1/4}}e(2\sqrt{nt} )\right|^4\d t\cr
&\ll_\e T^{-2}+T^{-1+\varepsilon }
\int_T^{2T}\left|\sum_{n\sim N}\frac{d(n)}{n^{1/4}}e(2\sqrt{nt} )\right|^4\d t\cr}
$$
for some $1\ll N\ll y.$ Therefore we have ($a\sim b$ means that $b \le a \le b' \le 2b$)
$$
\eqalign{
&\int_T^{2T}|\Delta_1^{\prime}(t,y)|^4\d t\ll_\e  T^{-2}+T^{-1+\varepsilon }
\int_T^{2T}\left|\sum_{n\sim N}\frac{d(n)}{n^{1/4}}e(2\sqrt{nt} )\right|^4\d t\cr
&\ll_\e  T^{-2}+T^{-1+\varepsilon }\sum_{n_1,n_2, n_3,n_4\sim N}
\frac{d(n_1)d(n_2)d(n_3)d(n_4)}{(n_1n_2 n_3n_4)^{1/4}}\cr
& \times
\left|\int_T^{2T}e\Bigl(2(\sqrt{n_1}+\sqrt{n_2}-\sqrt{n_3}-\sqrt{n_4}\,)  t\Bigr)\d t\right|\cr
&\ll_\e \frac{T^{-1+\varepsilon}}{N}\sum_{n_1,n_2, n_3,n_4\sim N}
\min\left(T,\frac{\sqrt T}{| \Omega|}\right),\cr}\leqno(5.13)
$$
Here we used trivial estimation and the first derivative test, and we set
$$
\Omega \;: =\;\sqrt{n_1}+\sqrt{n_2}-\sqrt{n_3}-\sqrt{n_4}.
$$

Note that $\min\bigl(T,\sqrt {T}/|\Omega|\bigr)= T$ if $|\Omega|\le T^{-1/2}$.
In this case the contribution to the last sum in (5.13) is, by (3.10) of Lemma 6,

$$
\eqalign{
&\ll_\e \frac{T^{-1+\varepsilon}}{N}T (T^{-1/2}N^{7/2}+N^2)
\ll_\e (T^{-1/2}N^{5/2}+N)T^\varepsilon\cr
&\ll_\e (T^{-1/2}y^{5/2}+y)T^\varepsilon\ll_\e T^{-1/2+\varepsilon}y^{5/2},
\cr}
$$
on noting that $y\gg T^{1/2}.$  If $|\Omega| > T^{-1/2}$, then
$\min\bigl(T,\sqrt {T}/|\Omega|\bigr)= \sqrt {T}/|\Omega|$.
By Lemma 6 again, the contribution  is
$$
\eqalign{
&\ll_\e \max_{T^{-1/2}<\eta\ll N^{1/2}}\frac{T^{-1/2+\varepsilon}}{N\eta}
\sum_{|\Omega|\sim \eta}1\cr
&\ll_\e  \max_{T^{-1/2}<\eta\ll N^{1/2}}\frac{T^{-1/2+\varepsilon}}{N\eta}(\eta N^{7/2}+N^2)\cr
&\ll_\e (T^{-1/2}y^{5/2}+y)T^\varepsilon\ll_\e T^{-1/2+\varepsilon}y^{5/2}.
\cr}
$$

Inserting the above two estimates into (5.13) we obtain
$$
\int_T^{2T}|\Delta_1^{\prime}(t,y)|^4\d t \;\ll_\e\;   T^{-1/2+\varepsilon}y^{5/2}.
\leqno(5.14)
$$

Now we bound $\int_{12}^{*}.$ When $k=2,3,4,$ by H\"older's inequality,
(5.11), (5.12) and Lemma 1 we have
$$
\eqalign{
\int_{12}^{*}& = \int_T^{2T}\D^{k-1}(t,y)\D'(t,y)E(t)\d t\cr&
\ll \left(\int_T^{2T}|\Delta_1^{\prime}(t,y)|^2dt\right)^{\frac 12}
 \left(\int_T^{2T}|\Delta_1(t,y)|^{2k}\d t\right)^{  \frac{k-1}{2k}}\cr&
 \times \left(\int_T^{2T}|E(t)|^{2k}\d t\right)^{  \frac{1}{2k}}
\ll_\e T^{\frac k4+\frac 34+\varepsilon}y^{\frac 14} .\cr}\leqno(5.15)
$$

When $k=5,6,7,8,$ by H\"older's inequality again, (5.11), (5.14) and Lemma 1 we have
$$
\eqalign{
\int_{12}^{*}&\ll \left(\int_T^{2T}|\Delta_1^{\prime}(t,y)|^4\d t\right)^{\frac 14}
\times \left(\int_T^{2T}|\Delta_1(t,y)|^{\frac{4k}{3}}\d t\right)^{  \frac{3k-3}{4k}}\cr&
 \times \left(\int_T^{2T}|E(t)|^{\frac{4k}{3}}\d t\right)^{  \frac{3}{4k}}
\ll_\e T^{\frac 58+\frac 34M(\frac{4k}{3})+\varepsilon}y^{\frac 58}.\cr}\leqno(5.16)
$$
In (3.4) we have $M(A) = A/4$ for $A\le 262/27 = 9,\overline{703}$, and 
$M(A)= 131(A-2)/416$ for $262/27 \le A \le 11$. Thus by Lemma 2, 
inserting (5.15) and (5.16) into (5.9) we obtain
$$
 \int_{12}\;\ll_\e\; T^{(k+1)\theta+\varepsilon}+\left\{\aligned
  T^{\frac k4+\frac 34+\varepsilon}y^{\frac 14},&\;\roman{if}\; k=2,3,4,\\
 T^{\frac 58+\frac k4+\varepsilon}y^{\frac 58},&\; \roman{if} \;k=5,6,7, \\
 T^{\frac {171}{64}+\varepsilon}y^{\frac 58},& \;\roman{if}\; k=8. \endaligned
 \right.\leqno(5.17)
$$

\medskip
Now we evaluate $\int_{11}$ (see (5.8)). Using $\Delta_1(t,y)=\Delta(t)-\Delta_2(t,y),$ we have
$$
\int_{11}= \int_T^{2T}\D_1^k(t,y)(\log t+C)\d t = \int_4 + \,O\Biggl\{\int_5+\int_6\Biggr\},
$$
say, where
$$
\eqalign{
&\int_4=\int_T^{2T}\Delta^{k}(t)(\log t+C)\d t,\cr
&\int_5=\int_T^{2T}|\Delta^{k-1}(t)\Delta_2(t,y)|(\log t+C)\d t,\cr
&\int_6=\int_T^{2T}|\Delta_2^{k}(t,y)|(\log t+C)\d t.
\cr}\leqno(5.18)
$$

From (5.4) and (5.6) we infer that
$$
\int_6\ll_\e {\left(\frac{T^{1/2+\varepsilon}}{y^{1/2}}\right)}^{k-2}\int^{2T}_T|\Delta_2 (t,y)|^2\d t
\ll_\e T^{\frac{k+1}{2}+\varepsilon}y^{-\frac{k-1}{2}}.
$$
By Cauchy's inequality, (5.6) and Lemma 2 we have, if $k=2,3,4,5$, that
$$
\eqalign{
\int_5&\ll \log T \left(\int_T^{2T}|\Delta_2 (t,y)|^2\d t\right)^{1/2}
\left(\int_T^{2T}|\Delta (t)|^{2k-2}dt\right)^{1/2}\cr
&\ll_\e T^{1+k/4+\varepsilon}y^{-1/4}.
\cr}\leqno(5.19)
$$
Similarly we obtain by H\"older's inequality, when $k=6,7,8$,
$$
\eqalign{
\int_5&\ll \log T \left(\int_T^{2T}|\Delta_2 (t,y)|^4\d t\right)^{1/4}\left(\int_T^{2T}
|\Delta (t)|^{4(k-1)/3}\d t\right)^{3/4}\cr
&\ll_\e  T^\varepsilon \left(\frac{T}{y}\int_T^{2T}|\Delta_2 (t,y)|^2\d t\right)^{1/4}
\left(\int_T^{2T}|\Delta (t)|^{4(k-1)/3}\d t\right)^{3/4}\cr
&\ll_\e T^\varepsilon \left(\frac{T^{5/2}}{y^{3/2} } \right)^{1/4}\left(T^{1+(k-1)/3}\right)^{3/4}
\cr
&= T^{9/8+k/4+\varepsilon}y^{-3/8},
\cr}\leqno(5.20)
$$
where we used $M(4(k-1)/3) = (k-1)/3$ by Lemma 2, since $4(k-1)/3 \le 28/3$.
Namely, for $k\le 8$ we have $(4k-4)/3 \le 28/3$, and by (3.4) with $\t = 131/416$ one
obtains $M(A) = A/4$ for $A \le 262/27 = 9.70370\ldots,$ while $28/3 = 9.3333\ldots$.

\medskip
Combining (5.19) and (5.20) with the above estimate for $\int_6$, we obtain
$$
 \int_{5}+\int_6 \;\ll_\e\;  \left\{\aligned
T^{1+k/4+\varepsilon}y^{-1/4},&  \;\roman{if}\; k=2,3,4,5,\\
T^{9/8+k/4+\varepsilon}y^{-3/8},&  \;\roman{if} \; k=6,7,8.
\endaligned\right.\leqno(5.21)
$$

From (5.1), (5.17) and (5.21) we have
$$
\int_1=\int_T^{2T}\Delta^{k}(t)(\log t+C)\d t+
O_\e\Bigl\{G_{k1}(T,y)T^\varepsilon+G_{k2}(T,y)T^\varepsilon\Bigr\},\leqno(5.22)
$$
say, where we have set
$$
 G_{k1}(T,y):=  \left\{\aligned
T^{1+k/4+\varepsilon}y^{-1/4},&\;\roman{if}\; k=2,3,4,5,\\
T^{9/8+k/4+\varepsilon}y^{-3/8},& \;\roman{if} \;k=6,7,8,
\endaligned\right.\leqno(5.23)
$$
and
$$
G_{k2}(T,y):= T^{(k+1)\theta+\varepsilon}+\left\{\aligned
  T^{\frac k4+\frac 34+\varepsilon}y^{\frac 14},&\;\roman{if} \;k=2,3,4,\\
 T^{\frac 58+\frac k4+\varepsilon}y^{\frac 58},& \;\roman{if}\; k=5,6,7, \\
  T^{\frac {171}{64}+\varepsilon}y^{\frac 58},& \;\roman{if}\; k=8.
\endaligned\right.\leqno(5.24)
$$

Now we estimate $\int_2$ (see (5.5)). Taking $k=2,4,6$ in the estimate  (5.22)--(5.24) we obtain
$$
\int_T^{2T}|\Delta_1(t,y)|^{k} |\zt|^2 \d t
\ll_\e T^{1+k/4+\varepsilon}\quad( k=2,4,6; \; T^{1/2}\ll  y\ll T^{3/5}).\leqno(5.25)
$$
Similarly, taking $k=8$ in (5.22)--(5.24) we obtain
$$
\int_T^{2T}|\Delta_1(t,y)|^8 |\zt|^2\d t\ll_\e T^{3+\varepsilon}\quad(T^{1/2}\ll  y\ll T^{21/40}),
\leqno(5.26)
$$
which combined with
  H\"older's inequality implies, for any
$2\le A\le 8$, that
$$
\eqalign{
&\int_T^{2T}|\Delta_1(t,y)|^A |\zt|^2 \d t \cr
&\ll \left( \int_T^{2T}|\Delta_1(t,y)|^8 |\zt|^2 \d t\right)^{\frac A8}
\left(\int_T^{2T}|\zt|^2\d t\right)^{1-\frac A8}  \cr
&\ll_\e T^{1+A/4+\varepsilon}\cr}\leqno(5.27)
$$
if $\; T^{1/2}\ll y\ll T^{21/40}.$

\medskip
When $k=2,3,4,$ from (5.7) with $k=2$, (5.25) and the Cauchy-Schwarz inequality we obtain, for
$T^{1/2}\ll  y\ll T^{3/5}$, that
$$
\eqalign{
\int_2&\ll \left(\int_T^{2T}|\Delta_2(t,y)|^{2}|\zt|^2\d t\right)^{1/2}\cr
&\times \left(\int_T^{2T}|\Delta_1(t,y)|^{2k-2} |\zt|^2\d t\right)^{1/2}\cr
&\ll_\e T^{\frac 98+\frac k4+\varepsilon}y^{-\frac 38}.\cr}
$$

\medskip
When $k=5,6,7,$ from (5.7) with $k=4$, (5.25) and H\"older's inequality we have,
for
$T^{1/2}\ll  y\ll T^{21/40}$, that
$$
\eqalign{
\int_2&\ll \left(\int_T^{2T}|\Delta_2(t,y)|^{4}|\zt|^2\d t\right)^{1/4}\cr
&\times \left(\int_T^{2T}|\Delta_1(t,y)|^{4(k-1)/3} |\zt|^2\d t\right)^{3/4}\cr
&\ll_\e T^{\frac{19}{16}+\frac k4+\varepsilon}y^{-\frac {7}{16}}.\cr}
$$

\medskip
When $k=8,$ from (5.7) with  with $k=8$, (5.27) with $A=8$ and H\"older's inequality we have, for
$T^{1/2}\ll  y\ll T^{21/40}$,
$$
\eqalign{
\int_2&\ll \left(\int_T^{2T}|\Delta_2(t,y)|^8|\zt|^2\d t\right)^{1/8}\cr
&\times \left(\int_T^{2T}|\Delta_1(t,y)|^{8} |\zt|^2\d t\right)^{7/8}\cr
&\ll_\e T^{\frac{103}{32}+ \varepsilon}y^{-\frac {15}{32}}.\cr}
$$

By combining the above three estimates it follows that
$$
\int_2\ll G_{k3}(T,y):=\left\{\aligned
 T^{\frac 98+\frac k4+\varepsilon}y^{-\frac 38},&\;\roman{when}\; k=2,3,4,\\
T^{\frac{19}{16}+\frac k4+\varepsilon}y^{-\frac {7}{16}} ,& \;\roman{when} \;k=5,6,7,\\
 T^{\frac{103}{32}+ \varepsilon}y^{-\frac {15}{32}} ,&\;\roman{when} \;k=8.
\endaligned\right.\leqno(5.28)
$$

From (5.5), (5.7), (5.22) and (5.28) we have
$$
\eqalign{
&\int_T^{2T}\Delta^{k}(t)|\zeta(\hf+it)|^2\d t =\int_T^{2T} \Delta^k(t)(\log t+C)\d t \cr
& +O_\e\left( \sum_{j=1}^3G_{kj}(T,y)T^\varepsilon+T^{\frac{2k+3}{4}+
\varepsilon}y^{-\frac{2k-1}{4}}\right),\cr}\leqno(5.29)
$$
where $G_{kj}(T,y)\;(j=1,2,3)$ was defined in (5.23), (5.24) and (5.28), respectively.
It is easy to see that
$$
\sum_{j=1}^3G_{kj}(T,y)\ll  \left\{\aligned
T^{\frac{9}{8}+\frac{k}{4}}y^{-\frac{3}{8}}+T^{\frac{3}{4}+\frac{k}{4}}y^{ \frac{1}{4}}+T^{(k+1)\theta},&
\;\roman{when}\; k=2,3,4,\\
T^{\frac{19}{16}+\frac{k}{4}}y^{-\frac{7}{16}}+T^{\frac{5}{8}+
\frac{k}{4}}y^{ \frac{5}{8}}+T^{(k+1)\theta},&\;\roman{when} \;k=5,6,7, \\
T^{\frac{103}{32} }y^{-\frac{15}{32}}+T^{\frac{171}{64} }y^{ \frac{5}{8}}+T^{9\theta},
&\;\roman{when} \;k=8.
\endaligned\right.\leqno(5.30)
$$
Now taking
$$
y=  \left\{\aligned
T^{\frac{3}{5} } ,&\;\roman{when}\; k=2,3,4,\\
T^{\frac{21 }{40} },& \;\roman{when} \;k=5,6,7,\\
T^{\frac{1}{2} },&\;\roman{when} \;k=8,
\endaligned\right.
$$
we obtain
$$
 \sum_{j=1}^3G_{kj}(T,y) +T^{\frac{2k+3}{4} }y^{-\frac{2k-1}{4}}\ll T^{1+\frac k4-\eta_k^{**}},
 \leqno(5.31)
$$
where
$$
 \eta_k^{**}: =  \left\{\aligned
1/10,&\;\roman{when} \;k=2,3,4,\\
27/640,&\;\roman{when} \;k=5,6,7, \\
1/64 ,&\;\roman{when} \;k=8.
\endaligned\right.\leqno(5.32)
$$
In the case when $k=2,3,4,8$ we equalize the terms containing $y$ in (5.30),
and $(k+1)\t < 1 + k/4 - \eta_k^{**}$ holds. In the case when $k = 5,6,7,$ note
that $T^{\frac{19}{16}+\frac{k}{4}}y^{-\frac{7}{16}} \ge T^{\frac{5}{8}+
\frac{k}{4}}y^{ \frac{5}{8}}$ for $T^{1/2} \le y \le T^{9/17}$ but as
$y \ll T^{21/40}$ has to hold and $21/40 < 9/17$, we take $y= T^{\frac{21 }{40}}$
to obtain (5.31) in this case as well.
\medskip
From (5.29)-(5.32) we obtain
$$
  \int_T^{2T}\Delta^{k}(t)|\zeta(\hf+it)|^2\d t =
  \int_T^{2T} \Delta^k(t)(\log t+C)\d t   +O_\e\left(T^{1+\frac k4-\eta_k^{**}+\varepsilon} \right),
$$
which  implies that
$$
  \int_1^{T}\Delta^{k}(t)|\zeta(\hf+it)|^2\d t =
  \int_1^{T} \Delta^k(t)(\log t+C)dt   +O_\e\left(T^{1+\frac k4-\eta_k^{**}+\varepsilon} \right).
  \leqno(5.33)
$$
From (5.33),  (2.1) and integration by parts we have ($\eta^*_k \equiv c(k)$)
$$
\eqalign{
&\int_1^{T}\Delta^{k}(t)|\zeta(\hf+it)|^2\d t\cr
& =C_k\Bigl(1+\frac k4\Bigr)\int_1^{T} t^{\frac k4}(\log t+C)\d t
   +O_\e\left(T^{1+\frac k4-\eta_k^{*}+\varepsilon}
 +T^{1+\frac k4-\eta_k^{**}+\varepsilon} \right)\cr
 &=C_kT^{1+\frac k4}\left(\log T+C-\frac{4}{k+4}\right) +O_\e\left(T^{1+\frac k4-\eta_k^{*}+\varepsilon}
 +T^{1+\frac k4-\eta_k^{**}+\varepsilon} \right) \cr
 &=c_1(k)T^{1+\frac k4}\log T + c_2(k)T^{1+\frac k4} +
 O_\e\left(T^{1+\frac k4-\eta_k +\varepsilon}\right),
\cr}
$$
where
$$
\eqalign{
c_1(k) \;&=\; C_k,\quad c_2(k) = C_k\left(C-\frac{4}{k+4}\right),\cr
\eta_k \;&=\; \min(\eta_k^{*},\,\eta_k^{**}) \quad (2\le k\le 8),\cr}
$$
so that
$$
\eta_2=  \eta_3=  \eta_4=1/10,\ \eta_5=3/80,\ \eta_6=35/4742,\ \eta_7=17/6312,\ \eta_8=8/9433.
$$
This ends the proof of Theorem 2.

\medskip
\head
6. Proof of Theorem 3
\endhead

\medskip
We retain the notation of Section 5. The main task is to evaluate
$$
\int_{11} := \int_T^{2T}\Delta_1^{2}(t,y)(\log t+C)\d t
$$
and to bound
$$
\int_{12}^{*}:=\int_T^{2T}\Delta_1(t,y)\Delta_1^{\prime}(t,y)E(t)\d t.
$$
By using (5.3) we have
$$
\eqalign{
\int_{11} &= \frac{1}{2\pi^2}\sum_{m,n\le y}\frac{d(m)d(n)}{(mn)^{3/4}}\cr&
\times\int_T^{2T}t^{1/2}(\log t+C)\cos(4\pi\sqrt{mt}-\pi/4)\cos(4\pi\sqrt{nt}-\pi/4)\d t.
\cr}\leqno(6.1)
$$
We use the identity
$$
\cos\a\cos\b = \hf\Bigr(\cos(\a+\b) + \cos(\a-\b)\Bigl)
$$
with $\a = 4\pi\sqrt{mt}-\pi/4, \b = 4\pi\sqrt{nt}-\pi/4$. The terms coming
from $\cos(\a+\b)$ make, by the first derivative test, a contribution which
is $\ll T\log^5T$. The same bound holds for the terms coming from
$\cos(\a-\b)$ when $m\ne n$. Finally, the terms $m=n$ contribute
$$
\eqalign{&
\frac{1}{4\pi^2}\sum_{n\le y}d^2(n)n^{-3/2}\int_T^{2T}t^{1/2}(\log t+C)\d t\cr&
= \frac{1}{4\pi^2}\sum_{n=1}^\infty d^2(n)n^{-3/2}\int_T^{2T}t^{1/2}(\log t+C)\d t
+ O(T^{3/2}y^{-1/2}\log^4T).\cr}
$$
It follows that
$$
\eqalign{
\int_{11} &= \frac{1}{4\pi^2}\sum_{n=1}^\infty d^2(n)n^{-3/2}\int_T^{2T}t^{1/2}(\log t+C)\d t
\cr&
+ O(T^{3/2}y^{-1/2}\log^4T) + O(T\log^5T).
\cr}\leqno(6.2)
$$
Now we estimate $\int_{12}^{*}$. Here we use the method of proof of Theorem 1 and replace
$E(t)$ by Lemma 9 (Atkinson's formula with $N=T$). We write
$$
\eqalign{
\int_{12}^{*} &= \int_{121}^{*} + \int_{122}^{*},\cr
\int_{121}^{*} &= \int_T^{2T}\D_1(t,y)\D_1'(t,y)\sum\nolimits_1(t)\d t,\cr
\int_{122}^{*} &= \int_T^{2T}\D_1(t,y)\D_1'(t,y)
\Bigl(\sum\nolimits_2(t)+ O(\log^2t)\Bigr)\d t.
\cr}\leqno(6.3)
$$
By H\"older's inequality we obtain
$$
\eqalign{
\int_{122}^{*} & \ll {\Bigl\{\int_T^{2T}
\Bigl(\sum\nolimits_2(t)+ O(\log^2t)\Bigr)^2\d t\Bigr\}}^{1/2}\cr&
\times {\Bigl\{\int_T^{2T}
|\D_1'(t,y)|^3\d t\Bigr\}}^{1/3}{\Bigl\{\int_T^{2T}
|\D_1(t,y)|^6\d t\Bigr\}}^{1/6}
\cr&
\ll_\e T^{11/12+\e}y^{1/2}.
\cr}\leqno(6.4)
$$
Here we bounded the mean square of $\Sigma_2(t)$ as after (4.6), used the bound
$$
\int_T^{2T}|\D_1'(t,y)|^3\d t \ll_\e T^\e y^{3/2},
$$
which follows from the Cauchy-Schwarz inequality from (5.12) and (5.14),
and (5.11) with $A = 6, M(6) = 3/2$.

\medskip
Having in mind (5.10), we see that the major contribution to $\int_{121}^{*}$
comes from a multiple of
$$
\eqalign{&
\int_T^{2T}\sum\nolimits_1(t)\sum_{n\le y}d(n)n^{-1/4}
\sin(4\pi\sqrt{nt}-\textstyle{\frac{\pi}{4}})\cr&
\times
\sum\limits_{m\le y}d(m)m^{-3/4}
\cos(4\pi\sqrt{mt}-\textstyle{\frac{\pi}{4}})\d t.\cr}\leqno(6.5)
$$
We use the explicit expression for $\sum_1$ given by Lemma 9.
By a splitting argument one sees that the integral in (6.5) can be written as $O(\log^3T)$
integrals of the form
$$
\eqalign{
I &:= I(T;L,M,N) =
\int_T^{2T}\frac{t^{1/4}}{L^{3/4}M^{3/4}N^{1/4}}\sum_{L<\ell\le2L}
c_1(\ell)e(t,\ell)\cos(f(t,\ell)\cr&
\times
\sum_{M<m\le2M}c_2(m)\cos(4\pi\sqrt{mt}-\pi/4)
\sum\limits_{N<n\le2N}c_3(n)\sin(4\pi\sqrt{nt}-\pi/4)\d t,
\cr}
$$
say, where the coefficients $c_j$ satisfy
$$
c_1(\ell) \ll d(\ell),\quad c_2(m) \ll d(m),\quad c_3(n) \ll d(n),
$$
and the functions $e(t,\ell), f(t,\ell)$ are as in Lemma 9.

We consider separately several cases.

{\bf Case 1}. $L \ge 100\max(M,N)$.

\medskip
In this case $I$ can be written as a linear combination of integrals
$$
\eqalign{
I' &= \frac{1}{L^{3/4}M^{3/4}N^{1/4}}\sum_{L<\ell\le2L}c_1(\ell)\sum_{M<m\le2M}c_2(m)
\sum\limits_{N<n\le2N}c_3(n)\cr&
\times
\int_T^{2T}t^{1/4}e(t,\ell)\exp\Bigl(if(t,\ell) \pm4\pi i\sqrt{mt}\pm4\pi i\sqrt{nt}\Bigr)
\d t.\cr}
$$
Then the derivative of the function in the exponential is $\gg \sqrt{\ell/T}$ and
(this is similar to the discussion regarding (4.8) and (4.9)), by the first derivative test,
we obtain
$$
I \;\ll\; T^{3/4}y^{3/4}\log^3T.\leqno(6.6)
$$

{\bf Case 2}. $M \ge 100\max(L,N)$.

\medskip

{\bf Case 3}. $N \ge 100\max(L,M)$.

\medskip
These cases are analogous to Case 1, and thus the analogue of (6.6) will hold.

\medskip
{\bf Case 4}. $N < 100\max(L,M),\; L < 100\max(N,M),\; M < 100\max(L,N).$

\medskip In this case, like in (4.10) in the proof of Theorem 1, we shall use
mean value estimates. To this end let
$$
\eqalign{
U_1(t) &:= \sum_{L<\ell\le2L}c_1(\ell)e(t,\ell)\cos\bigl(f(t,\ell)\bigr),
\cr
U_2(t)&:= \sum_{M<m\le2M}c_2(m)\cos(4\pi\sqrt{mt}-\pi/4),\cr
U_3(t) &:= \sum\limits_{N<n\le2N}c_3(n)\sin(4\pi\sqrt{nt}-\pi/4).
\cr}
$$
We need the bounds
$$
\int_T^{2T}|U_1(t)|^2\d t \;\ll\;TL\log^4T\leqno(6.7)
$$
and
$$
\int_T^{2T}|U_1(t)|^4\d t \;\ll_\e\;T^\e(T^{1/2}L^{7/2} + TL^2).\leqno(6.8)
$$
Note that (6.7) follows directly by squaring out the integrand and integrating,
while (6.8) follows by the use of (3.10) of Lemma 6 with $k =2$,
similarly as in (5.13) in the proof of Theorem 2. We also note that the analogues
of (6.7) and (6.8) hold for the corresponding integrals of $U_j(t)$ ($j = 2,3$).

\medskip If Case 4 holds, then we must have

\medskip {\bf Case 4.1}. $L \ll M, L\ll N, M\asymp N$, or

\medskip {\bf Case 4.2}. $M \ll L, M \ll N, L\asymp N$, or

\medskip {\bf Case 4.3}. $N \ll L, N \ll M, L\asymp M$.

\medskip
Let us consider first the case 4.1. Using (6.7)--(6.8) and its analogues, and
H\"older's inequality, we have
$$
\eqalign{
I& \ll \frac{T^{1/4}}{L^{3/4}M^{3/4}N^{1/4}}\int_T^{2T}|U_1(t)U_2(t)U_3(t)|\d t\cr&
\ll \frac{T^{1/4}}{L^{3/4}M^{3/4}N^{1/4}}{\Bigl(\int_T^{2T}|U_1(t)|^4\d t\Bigr)}^{\frac14}
{\Bigl(\int_T^{2T}|U_2(t)|^4\d t\Bigr)}^{\frac14}
{\Bigl(\int_T^{2T}|U_3(t)|^2\d t\Bigr)}^{\frac12}\cr&
\ll_\e \frac{T^{1/4+\e}}{L^{3/4}M^{3/4}N^{1/4}}
(T^{1/8}L^{7/8}+T^{1/4}L^{1/2})(T^{1/8}M^{7/8}+T^{1/4}M^{1/2})T^{1/2}N^{1/2}\cr&
\ll_\e T^\e(TL^{1/8}M^{3/8} + T^{9/8}L^{1/8} + T^{9/8}M^{3/8}L^{-1/4}+ T^{5/4}L^{-1/4}).
\cr}
$$
By using the trivial estimate $U_1(t) \ll L\log L$ we also have, since $M\asymp N$,
$$
\eqalign{
I& \ll \frac{T^{1/4}L^{1/4}\log L}{M}\int_T^{2T}|U_2(t)U_3(t)|\d t\cr&
\ll \frac{T^{1/4}L^{1/4}\log L}{M}{\left(\int_T^{2T}|U_2(t)|^2\d t
\int_T^{2T}|U_3(t)|^2\d t\right)}^{1/2}\cr&
\ll T^{5/4}L^{1/4}\log^5T.
\cr}
$$
From the last two estimates for $I$ we infer that, when $L\ll M, L\ll N, M\asymp N$,
$$
\eqalign{
I& \ll_\e T^\e \left(TL^{1/8}M^{3/8} + T^{9/8}L^{1/8} + T^{5/4}L^{-1/4}
+ \min\Bigl(\frac{T^{9/8}M^{3/8}}{L^{1/4}},\,T^{5/4}L^{1/4}\Bigr)\right)\cr&
\ll_\e T^\e(Ty^{1/2} + T^{19/16}y^{3/16} + T^{5/4}),
\cr}\leqno(6.9)
$$
since $\min(a,b) \le \sqrt{ab}$ for $a,b>0$.

\medskip In the case 4.2, the argument is the same, only the orders of $L$  and $M$ are changed.
Consequently the bound (6.9) will hold again. Finally in the case 4.3 we obtain

$$
\eqalign{
I& \ll \frac{T^{1/4}}{L^{3/4}M^{3/4}N^{1/4}}
{\Bigl(\int_T^{2T}|U_1(t)|^4\d t\Bigr)}^{1/4}
\cr&\times
{\Bigl(\int_T^{2T}|U_2(t)|^4\d t\Bigr)}^{1/4}
{\Bigl(\int_T^{2T}|U_3(t)|^2\d t\Bigr)}^{1/2}\cr&
\ll_\e \frac{T^{1/4+\e}}{L^{3/4}M^{3/4}N^{1/4}}
T^{1/2}N^{1/2}(T^{1/8}M^{7/8}+T^{1/4}M^{1/2})(T^{1/8}L^{7/8}+T^{1/4}L^{1/2})\cr&
\ll_\e \frac{T^{3/4+\e}N^{1/4}}{M^{3/2}}(T^{1/4}M^{7/4}+T^{1/2}M)\cr&
\ll_\e T^{1+\e}N^{1/4}M^{1/4} + T^{5/4+\e}N^{1/4}M^{-1/2} \ll_\e T^{1+\e}y^{1/2}+T^{5/4+\e}.
\cr}
\leqno(6.10)
$$
Hence (6.9) and (6.10) yield
$$
I \ll_\e T^\e(Ty^{1/2} + T^{19/16}y^{3/16} + T^{5/4}).
$$
Combining the estimates for $I$ in all four cases we have
$$
I \ll_\e T^\e(T^{3/4}y^{3/4}+Ty^{1/2} + T^{19/16}y^{3/16} + T^{5/4})
\ll_\e T^\e(Ty^{1/2} + T^{19/16}y^{3/16} + T^{5/4}),
\leqno(6.11)
$$
since $T^{3/4}y^{3/4} \le Ty^{1/2}$. Using (6.11) to 
bound the expression in (6.5) we obtain
$$
\int^*_{121} \;\ll_\e\; T^{1+\e}y^{1/2} + T^{19/16+\e}y^{3/16}+T^{5/4+\e}.\leqno(6.12)
$$
From   (5.8),  (6.2) and (6.12) we have
$$
\int_{12}^*  \;\ll_\e\; T^{1+\e}y^{1/2} + T^{19/16+\e}y^{3/16} + T^{5/4+\e},
\leqno(6.13)
$$
and this gives in (5.9) (note that $k=2$)
$$
\int_{12} \;\ll_\e\; T^{1+\e}y^{1/2} + T^{19/16+\e}y^{3/16} + T^{5/4+\e}.
\leqno(6.14)
$$
From (5.8), (6.2) and (6.14) it follows that
$$
\eqalign{
\int_1 &=
 \frac{1}{4\pi^2}\sum_{n=1}^\infty d^2(n)n^{-3/2}\int_T^{2T}t^{1/2}(\log t+C)\d t
\cr&
+ O_\e(T^{3/2+\e}y^{-1/2}) + O_\e(T^{1+\e}y^{1/2}) + O_\e(T^{19/16+\e}y^{3/16})
 + O_\e(T^{5/4+\e}).
\cr}\leqno(6.15)
$$
It remains yet to deal with $\int_2$ and $\int_3$ in (5.5) when $k=2$. In this case
(5.7) yields
$$
\int_3\;\ll_\e\;T^{7/4+\e}y^{-3/4}.\leqno(6.16)
$$
Now we bound $\int_2$. We have
$$
\int_T^{2T}|\D_2(t,y)|^3\d t \ll_\e T^{1/2+\e}y^{-1/2}\int_T^{2T}|\D_2(t,y)|^2\d t
\ll_\e T^{2+\e}y^{-1}.
$$
We use (5.11) with $A=6$, Lemma 3 and H\"older's inequality to obtain that
$$
\eqalign{
\int_2 &\ll {\Bigl(\int_T^{2T}|\zt|^4\d t\Bigr)}^{1/2}
 {\Bigl(\int_T^{2T}|\D_2(t,y)|^3\d t\Bigr)}^{1/3}
 \cr&
\times  {\Bigl(\int_T^{2T}|\D_1(t,y)|^6\d t\Bigr)}^{1/6}\cr&
  \ll_\e T^{19/12+\e}y^{-1/3}.\cr}
  \leqno(6.17)
  $$
Thus from  (6.15)--(6.17) it follows that
$$
\eqalign{&
\int_T^{2T}\D^2(t)|\zt|^2\d t
=  \frac{1}{4\pi^2}\sum_{n=1}^\infty d^2(n)n^{-3/2}\int_T^{2T}t^{1/2}(\log t+C)\d t
\cr&
+ O_\e\left(T^\e(T^{19/12}y^{-1/3}+  Ty^{1/2} + T^{19/16}y^{3/16} + T^{5/4})\right),
\cr}\leqno(6.18)
$$
since $T^{3/2}y^{-1/2} \ll Ty^{1/2}$ for $y \gg T^{1/2}$. Finally, taking
$$
y \;=\; T^{7/10}
$$
it is seen that all the error terms in (6.18) are $\ll_\e T^{27/20+\e}$,
and we obtain from (6.18)
$$
\eqalign{&
\int_1^{T}\D^2(t)|\zt|^2\d t \cr&
=  \frac{1}{4\pi^2}\sum_{n=1}^\infty d^2(n)n^{-3/2}\int_1^{T}t^{1/2}(\log t+C)\d t
+ O_\e(T^{27/20+\e})
\cr&
= c_1(2)T^{3/2}\log T + c_2(2)T^{3/2} + O_\e(T^{3/2-3/20+\e}),
\cr}
$$
which is the assertion of Theorem 3.

\medskip
\Refs
\medskip

\item{[1]} F.V. Atkinson, {\it The mean value of the Riemann zeta-function},
Acta Math. {\bf81}(1949), 353-376.

\item{[2]} D.R. Heath-Brown,
{\it The distribution and moments of the error term in the Dirichlet
divisor problems}, Acta Arith. {\bf60}(1992), 389-415.

\item  {[3]} M. N. Huxley, {\it Exponential sums and lattice points III,}
Proc. London Math. Soc. {\bf87}(3)(2003), 591--609.

\item{[4]} A. Ivi\'c, {\it The Riemann zeta-function}, John Wiley \&
Sons, New York, 1985 (2nd ed. Dover, Mineola, New York, 2003).

\item{[5]} A. Ivi\'c, {\it On some integrals involving the mean square
formula for the
    Riemann zeta-function}, Publs. Inst. Math. (Belgrade)
    {\bf46(60)} (1989), 33-42.

\item{[6]} A. Ivi\'c, {\it Large values of certain number-theoretic
error terms}, Acta  Arith.
  {\bf56} (1990), 135-159.

\item {[7]} A. Ivi\'c,  {\it Mean values of the Riemann zeta-function},
LN's {\bf 82},  Tata Inst. of Fundamental Research,
Bombay,  1991 (distr. by Springer Verlag, Berlin etc.).

\item {[8]} A. Ivi\'c, {\it On the divisor function and the Riemann 
zeta-function in short intervals},
The Ramanujan Journal, {\bf19}(2009), 207-224.

\item{[9]} A. Ivi\'c, {\it On some mean value results for the zeta-function
and a divisor problem}, to appear in Filomat, preprint available
at {\tt arXiv:1406.0604}.

\item{[10]} A. Ivi\'c and Y. Motohashi, {\it On the fourth power moment of the
Riemann zeta-function}, J. Number Theory {\bf51}(1995), 16-45.

\item{[11]} A. Ivi\'c and P. Sargos, {\it On the higher power moments of the
 error term in the divisor problem}, Illinois J. Math. {\bf81}(2007), 353-377.

\item{[12]} M. Jutila, {\it Riemann's zeta-function and the divisor problem},
Arkiv Mat. {\bf21}(1983), 75-96 and II, ibid. {\bf31}(1993), 61-70.

\item{[13]} M. Jutila, {\it On a formula of Atkinson},
Topics in classical number theory, Colloq. Budapest 1981,
Vol. I, Colloq. Math. Soc. J\'anos Bolyai {\bf34}(1984), 807-823.

\item{[14]} K.-L. Kong, {\it Some mean value theorems for certain error
terms in analytic number theory}, Master's Thesis, University of Hong
Kong, Hong Kong, 2014, 64pp.

\item{[15]} Y.-K. Lau, K.-M. Tsang,
{\it On the mean square formula of the error term in the Dirichlet divisor problem},
Math. Proc. Cambridge Philos. Soc. {\bf146}(2009), no. 2, 277-287.

\item{[16]} T. Meurman, {\it On the mean square of the Riemann zeta-function},
Quart. J. Math. Ser. (2){\bf38}\break(1987), 337-343.

\item{[17]} Y. Motohashi, {\it Spectral theory of the Riemann zeta-function}, Cambridge
University Press, Cambridge, 1997.

\item{[18]} K.  Ramachandra, {\it On the mean-value and omega-theorems
for the Riemann zeta-function}, LN's {\bf85}, Tata Inst. of Fundamental Research
(distr. by Springer Verlag, Berlin etc.), Bombay, 1995.

\item{[19]} O. Robert and P. Sargos, {\it Three-dimensional exponential
sums with monomials, J. reine angew. Math.} {\bf591}(2006), 1-20.

\item{[20]} P. Shiu, {\it A Brun-Titchmarsh theorem for multiplicative functions},
J. reine angew. Math. {\bf31}(1980), 161-170.

\item{[21]} K.-M. Tsang, {\it Recent progress on the Dirichlet divisor problem and the
mean square of the Riemann zeta-function}, Sci. China Math. {\bf53}(2010), no. 9, 2561-2572.

\item{[22]} G.F. Vorono{\"\i}, Sur une fonction transcendante et ses
applications \`a la sommation de quelques s\'eries, Ann. \'Ecole Normale
{\bf21}(3)(1904), 207-268 and ibid. {\bf21}(3)(1904), 459-534.

\item{[23]} N. Watt, {\it A note on the mean square of $|\zt|$},
J. London Math. Soc. {\bf82(2)}(2010), 279-294.

\item{[24]} W. Zhai, {\it On higher-power moments of $\Delta(x)$}, Acta Arith.
{\bf112}(2004), 367-395; II. ibid. {\bf114}(2004), 35-54; III ibid. {\bf118}(2005),
263-281, and IV,
Acta Math. Sinica, Chin. Ser. {\bf49}\break (2006), 639-646.

\item{[25]} W. Zhai, {\it On higher-power moments of $E(t)$}, Acta Arith. {\bf115}(2004),
329-348.


\endRefs

\enddocument

\bye